\documentclass[12pt,reqno]{article}
\usepackage{amssymb}
\usepackage{amsmath}
\usepackage{amsthm}
\usepackage{amsfonts}
\usepackage{url}
\usepackage{breakurl}
\usepackage{hyperref}

\renewcommand{\ge}{\geqslant}
\renewcommand{\le}{\leqslant}

\newcommand{\C}{{\mathbb C}}

\newcommand{\circulant}{{\rm circ}}
\newcommand{\Sz}{S_{01}}
\newcommand{\Spm}{S_{\pm1}}
\newcommand{\HBE}{\mathrm{H_{BE}}}

\newcommand{\x}{{\varepsilon}}

\newcommand{\nsp}{\hspace{0em}}

\theoremstyle{plain}
\newtheorem{theorem}{Theorem}

\theoremstyle{definition}

\newcommand{\seqnum}[1]{\href{http://oeis.org/#1}{\underline{#1}}}

\begin{document}
\bibliographystyle{plain}
\title{Computation of Maximal Determinants of Binary Circulant Matrices}
\author{Richard P. Brent%
\footnote{
Mathematical Sciences Institute, Australian National University,
Canberra, ACT 2600, Australia,
and CARMA, University of Newcastle, Callaghan, NSW 2308, Australia.
\href{mailto:circulants@rpbrent.com}{\tt circulants@rpbrent.com}
}
\and
Adam B.\ Yedidia%
\footnote{Department of Electrical Engineering and Computer Science,
Massachusetts Institute of Technology,
Cambridge, Massachusetts, USA.
\href{mailto:adamyedidia@gmail.com}{\tt adamyedidia@gmail.com}
}
}

\date{}

\maketitle

\begin{abstract}
We describe algorithms for computing maximal determinants of binary
circulant matrices of small orders. 	
Here ``binary matrix'' means a matrix whose elements
are drawn from $\{0,1\}$ or $\{-1,1\}$.
We describe efficient parallel algorithms for the search, using
Duval's algorithm for generation of necklaces and the well-known
representation of the determinant of a circulant in terms of roots of unity.
Tables of maximal determinants are given for orders $\le 53$.
Our computations extend earlier results and disprove two plausible
conjectures.
\end{abstract}

\pagebreak[3]
\section{Introduction}				\label{sec:intro}

A \emph{circulant} matrix $A = (a_{j,k})$ of order $n$ is an $n \times n$
matrix whose elements $a_{j,k}$ depend only on $(k-j) \bmod n$. 
Thus, an $n \times n$ circulant is a
matrix of the form $A = (a_{(k-j) \bmod n})_{\,0 \le j, k < n}$.
Circulants arise in various applications in signal processing and
combinatorics, and have a close connection with Fourier transforms.
The set of all circulants of order $n$ (with elements in some fixed
ring $R$) form a commutative algebra, since the sum and product of two
circulants is a circulant, and it is easy to see that multiplication
of circulants is commutative.

We write $\circulant(a_0,a_1,\ldots,a_{n-1})$ for the circulant 
$(a_{(k-j) \bmod n})_{\,0\le j,k<n}$ whose 
first row is $(a_0,a_1,\ldots,a_{n-1})$.

By a \emph{binary} matrix we mean a matrix whose elements are in one of the
sets $\Sz := \{0,1\}$ or $\Spm := \{-1,1\}$. It will be clear from
the context which of these two cases is being considered.
A \emph{binary circulant} is a circulant matrix whose
elements are in $\Sz$ or $\Spm$.

There is a natural one-to-one correspondence between the integers\linebreak
\hbox{$\{0,1,\ldots 2^n-1\}$} and the binary circulant matrices of order
$n$.  More precisely, if $N\in \{0,1,\ldots,2^n-1\}$ has the representation
\[N = \sum_{j=0}^{n-1} 2^{n-1-j\,}b_j,\]
so may be written in binary as $b_0 \ldots b_{n-1}$,
we associate $N$ with
$\circulant(a_0,\ldots,a_{n-1})$,
where $a_j = b_j$ in the case of $\Sz$, and $a_j = 2b_j-1$ in
the case of $\Spm$.

The \emph{maximal determinant problem} is concerned with the maximal value
of $|\det A|$ for an $n\times n$ binary matrix $A$. 
The \emph{Hadamard bound}~\cite{Hadamard93} states that, in the case of
binary matrices $A$ over $\{\pm1\}$, we have 
\begin{equation}
|\det A| \le n^{n/2}.					\label{eq:Hadamard1}
\end{equation}
Moreover, Hadamard's inequality is sharp for infinitely many $n$, for
example powers of two (Sylvester~\cite{Sylvester67}) or $n$ of the form
$q+1$ where $q$ is a prime power and 
$q \equiv 3 \bmod 4$ (Paley~\cite{Paley33}).

There is a well-known 
connection between the determinants of
$\{0,1\}$-matrices of order $n$ and
$\{\pm1\}$-matrices of order $n+1$.
This implies that an $(n+1)\times (n+1)$ $ \{\pm1\}$-matrix 
always has determinant divisible
by $ 2^{n}$.
See \cite{Neubauer97} 
for details.
We give an example with $n=3$, starting with an $n\times n$ binary
matrix $B$ and ending with an $(n+1)\times(n+1)$ $\{\pm1\}$-matrix $A$,
with $\det A = 2^n\det(B)$.
\[
B = 
\left(\begin{tabular}{ccc}
$1$&$0$&$1$\\
$1$&$1$&$0$\\
$0$&$1$&$1$
\end{tabular}\right)
\begin{tabular}{c}
\begin{tiny}{ double}\end{tiny}\\[-5pt]
$\longrightarrow$
\end{tabular}
\left(\begin{tabular}{rrr}
$2$&$0$&$2$\\
$2$&$2$&$0$\\
$0$&$2$&$2$
\end{tabular}\right)
\]

\[
\begin{tabular}{c}
\begin{tiny}{ border}\end{tiny}\\[-5pt]
$\longrightarrow$
\end{tabular}
\left(\begin{tabular}{rrrr}
1&1&1&1\\
$0$&$2$&$0$&$2$\\
$0$&$2$&$2$&$0$\\
$0$&$0$&$2$&$2$
\end{tabular}\right)
\begin{tabular}{c}
\begin{tiny}{ subtract}\end{tiny}\\[-5pt]
$\longrightarrow$\\[-7pt]
\begin{tiny}{ first row}\end{tiny}
\end{tabular}
\left(\begin{tabular}{rrrr}
$1$&$1$&$1$&$1$\\
$-1$&$1$&$-1$&$1$\\
$-1$&$1$&$1$&$-1$\\
$-1$&$-1$&$1$&$1$
\end{tabular}\right) = A.
\]
The doubling step is the
only step where the determinant changes, and there it is multiplied
by $2^n$.

Thus, Hadamard's bound~\eqref{eq:Hadamard1} gives the bound
\begin{equation}
|\det B| = 2^{-n}|\det A| \le 2^{-n}(n+1)^{(n+1)/2},	\label{eq:Hadamard2}
\end{equation}
which applies for all $\{0,1\}$-matrices $B$ of order $n$.
We shall refer to both \eqref{eq:Hadamard1} and \eqref{eq:Hadamard2} as
\emph{Hadamard's inequality}, since it will be clear from the context
which inequality is intended.\footnote{In fact, Hadamard 
in~\cite{Hadamard93} proved a more
general inequality than~\eqref{eq:Hadamard1}, and as far as we are aware he
never stated~\eqref{eq:Hadamard2} explicitly.
A simple proof of~\eqref{eq:Hadamard1} is given by Cameron~\cite{Cameron06}.
}

The mapping from $\{0,1\}$-matrices to $\{\pm1\}$-matrices  
is reversible if we are allowed to
normalise the first row and column of the $\{\pm1\}$-matrix by changing 
the signs of rows/columns as necessary. 

The transformation illustrated above (or its reverse) 
does \emph{not} preserve any
circulant structure. 

\emph{Hadamard matrices} are square matrices with entries in $\Spm$ and
mutually orthogonal rows. The order of a Hadamard matrix is $1$, $2$, or
a multiple of $4$.  It is not known whether a Hadamard matrix of order
$4k$ exists for every positive integer $k$ (this is the 
\emph{Hadamard conjecture}).

Various constructions for Hadamard matrices use circulant matrices. For
example, the first Paley construction~\cite{Paley33} uses a circulant matrix
of order~$p$, where $p$ is a prime, $p \equiv 3 \bmod 4$, to
construct a Hadamard matrix of order \hbox{$p+1$}.
(The Paley construction also works for prime powers, e.g.~$27 = 3^3$,
but does not involve circulants in such cases.)
Fletcher, Gysin and
Seberry~\cite{Fletcher01} use two circulants and a border of width
two to construct Hadamard matrices.  The Williamson
construction~\cite{Williamson44} requires four matrices $A, \ldots, D$ which
satisfy certain conditions, and for computational reasons these matrices
are usually taken to be circulants.

Circulant matrices also play an important role in noisy convolutional
Gaussian channels. Given a channel in which the output vector is given by
the convolution of the input vector with a chosen mask vector, in the
presence of additive Gaussian noise, the choice of mask that maximizes the
mutual information of the channel in high-SNR regimes is the first row of a
$\{0,1\}$-circulant
with near-flat Fourier spectrum,
and this circulant is often one with maximal or close to maximal 
determinant. This has important applications in X-ray and gamma ray
astronomy, optics, and computational
imaging~\cite{Ables68,Asif15,Dicke68,Fenimore78,Levin07,Yedidia17}.

It is well-known that the (unnormalised) eigenvectors of
$\circulant(a_0,\ldots,a_{n-1})$ are given by
$v_j = (1,\omega^j,\omega^{2j},\ldots,\omega^{(n-1)j})^T$,
$0 \le j < n$,
where $\omega$ is a primitive $n$-th root of unity. For example,
in $\C$ we can take $\omega := \exp(2\pi i/n)$.
It follows that the eigenvalues are
\begin{equation}
\lambda_j = a_0 + a_1\omega^j + \cdots + a_{n-1}\omega^{(n-1)j},
	\; 0 \le j < n,				\label{eq:lambda}
\end{equation}
and the determinant is
\begin{equation}
\prod_{j=0}^{n-1}\lambda_j = \prod_{j=0}^{n-1} f(\omega^j),
						\label{eq:det_product}
\end{equation}
where
\[f(z) := \sum_{k=0}^{n-1} a_k z^k.\]
The polynomial $f(z)$ is called the \emph{associated polynomial} of
the circulant.
Also, $f(z)$ is called a \emph{Littlewood polynomial} if the
coefficients $a_k \in \{\pm1\}$, and a 
\emph{Newman polynomial} if the $a_k \in \{0,1\}$ and $a_0 = 1$.

If $A = \circulant(a_0,\ldots,a_{n-1})$ is nonsingular, then
\eqref{eq:det_product} gives
\[
\frac{\log|\det A|}{n} = \frac{1}{n}\sum_{j=0}^{n-1}
	\log|f(e^{2\pi ij/n})|\,.
\]
This may be regarded as a discrete analogue of the
\emph{Mahler measure}~\cite{Smyth08}
\[
m(f) := \int_0^1 \log|f(e^{2\pi it})|\,{\rm d}t\,. 
\]
Using~\eqref{eq:det_product} to compute $\det A$ for a circulant matrix $A$
takes $O(n^2)$ arithmetic operations, whereas Gaussian elimination does not
take advantage of the circulant structure and takes of order $n^3$
operations. If we are considering binary matrices, whose determinants are
integers, it is necessary to perform the operations in
$\C$ to sufficient precision to obtain a result with absolute error less
than $1/2$, so that the correct result can be found by rounding 
to the nearest integer.
{From} the Hadamard bounds \eqref{eq:Hadamard1}--\eqref{eq:Hadamard2},
this means that we may have to work with of order $n\log n$ bits of precision.

\pagebreak[3]
To avoid the problem of rounding errors altogether, 
we can work over a finite field.
If $p$ is a prime such that $p \equiv 1 \bmod n$, and
$\rho$ is a primitive root (mod $p$), then\footnote{It is not
necessary to know a primitive root (mod $p$).  We can choose a random $a$,
compute $\omega = a^{(p-1)/n}$, and check if
$1, \omega, \omega^2,\ldots, \omega^{n-1}$
are distinct (mod $p$). If not, reject $\omega$ and repeat with another 
random $a$.  In this way we work in a (small) group of order
$n$, instead of a (large) group of order $p-1$, and there is no need
to factor $p-1$.
The expected number of iterations is $n/\phi(n) = O(\log\log n)$.} 
\[\omega = \rho^{(p-1)/n} \bmod p\]
is a primitive $n$-th root of unity in the finite field $F_p$, and we
can use~\eqref{eq:det_product} to compute $\det A \bmod p$.
If $U$ is an  upper bound on $|\det A|$,
and \hbox{$p \ge 2U+1$}, then the result mod~$p$ is sufficient to determine
$\det A$. 
Thus, if we use a Hadamard bound for $U$,
the prime $p$ should have of order $n\log n$ bits.
Alternatively, we could use several smaller primes with a sufficiently 
large product, and reconstruct the
result using the Chinese Remainder Theorem.\footnote{Tests indicate
that, at least for $n \le 50$, it is faster to use a single prime.
One reason for
this is that the value $\det A$ needs to be reconstructed for each
circulant $A$, so the cost of the reconstruction steps is not
negligible.}

\section{Lyndon words and necklaces}

The usual definition of a
\emph{Lyndon word} is a nonempty string that is strictly smaller in 
lexicographic order than all of its proper rotations.  
Thus, the first six Lyndon words over $\Sz$ are
$0$, $1$, $11$, $101$, $111$, and $1111$.
Lyndon words were introduced by Shirshov~\cite{Shirshov53} (who called them
``regular words'') and Lyndon~\cite{Lyndon54} (who called the ``standard
lexicographic sequences'').

Since we consider words of a fixed length $n$, it is convenient to use 
the concept of a (binary) necklace~\cite{Weisstein-necklace}.
We say that $w = w_0\ldots w_{n-1}$ is a
\emph{necklace of length $n$} if $w$ is not larger 
(in the lexicographic order) than any of
its rotations.  This corresponds to Duval's
``representative of a class of words of
length $n$''~\cite[(3) on pg.~258]{Duval88},
where two words are said to be
in the same class if one is a rotation of the other. 

For example, according to our definition,
the six necklaces of length $4$ over $\Sz$
are $0000$, $0001$, $0011$, $0101$, $0111$, and $1111$.
It can be seen that, if we strip off leading zeros, we obtain the first six
Lyndon words. Thus, the concepts of ``Lyndon word'' and ``necklace''
are closely related, and algorithms for one may often by modified
to apply to the other.

The number $K(n)$ of necklaces of length $n$ over a binary alphabet is
\begin{equation}
K(n) = \frac{1}{n}\sum_{d|n}2^{n/d}\phi(d) \simeq {2^n}/{n},
	\label{eq:Lyndon_kt}
\end{equation}
where $\phi$ is Euler's phi function.
$K(n)$ is tabulated in OEIS A000031~\cite{A000031}.

If $A$ is a circulant, then $|\det A|$ is 
invariant under rotations of the first row $(a_0,\ldots,a_{n-1})$.
Thus, when searching for circulants of order~$n$ 
with maximal determinants, it is
sufficient to consider circulants whose first row is a necklace
of length~$n$.
{From}~\eqref{eq:Lyndon_kt}, this saves a factor of approximately~$n$.

In our computations we use two nontrivial algorithms related to
Lyndon words/necklaces.
One is the algorithm of Booth~\cite{Booth80}, which determines in linear
time if a word $w = w_0\ldots w_{n-1}$ is in fact a necklace.\footnote{We
use a simplified version of Booth's algorithm since we do not 
need to know the rotation that would convert $w$ into a necklace.}
Booth's algorithm is closely related to the initial phase of the
Knuth, Morris and Pratt fast pattern-matching algorithm~\cite{KMP-algorithm}.

The other algorithm that we use is Duval's algorithm~\cite{Duval88} which,
given a necklace of length~$n$, returns the
next necklace (of length~$n$) in lexicographic order\footnote{Duval's
paper~\cite{Duval88} considers Lyndon words but, using 
\cite[comment~(3)~on~pg.~258]{Duval88},
we easily get a similar algorithm for necklaces.},
in amortised 
(i.e.\ average) constant time, see~\cite{Berstel94}. Using Duval's algorithm
we can cycle through all necklaces of length~$n$ in time
$O(2^n/n)$.

Other algorithms could be used.  For example, Shiloach~\cite{Shiloach81}
gives an algorithm that reduces the number of comparisons used by Booth's
algorithm. We used Booth's algorithm because it was sufficient for our
purposes, and simpler to implement than Shiloach's algorithm.
The overall complexity of our algorithms is dominated by the time
required to evaluate determinants using~\eqref{eq:det_product},
not by the time required to check or enumerate necklaces.

\section{Fast evaluation of circulant determinants}

Standard algorithms of linear algebra, such as Gaussian elimination,
require of order $n^3$ operations to evaluate the determinant of
an $n \times n$ matrix $A$.  Using formula~\eqref{eq:det_product}, this
can be reduced to order $n^2$ if $A$ is a circulant.  In fact, using
the fast Fourier transform (FFT), $O(n \log n)$ operations suffice.

However, in our application we can do even better. Because Duval's algorithm
takes constant time (on average), the number of symbols that are changed
as we go from one necklace to the next is $O(1)$ 
on average.\footnote{We find experimentally that the mean number
of symbols changed is 
$2+O(n/2^n)$ as $n \to \infty$. 
The limiting value $2$ is the same as the
mean number of bits changed when counting up in binary.}
Thus, each $\lambda_j$ value given by~\eqref{eq:lambda}
can be updated in $O(1)$ operations (on average), and the determinant,
given by~\eqref{eq:det_product}, can be updated with $O(n)$ operations
(on average).  Since there are $\simeq 2^n/n$ necklaces of
length $n$, the computation of all the relevant determinants
can be done with $O(2^n)$ operations.  The cost of precomputing a table
of powers $\omega^{jk}\; (0 \le j, k < n)$ is negligible.

Note that we used the term ``operations'' rather than ``time'', because the
arithmetic operations need to be performed using of order $n\log n$ bits
of precision, as noted above.  Thus, the overall complexity
is $O(2^n M(n\log n))$, where $M(N)$ is the time required to multiply $N$-bit
numbers.

In theory, a slightly better complexity can be attained by using several
small primes and reconstructing the result via the Chinese Remainder
Theorem.  However, the cost of $O(2^n/n)$ reconstructions must be taken into
account.  In practice, $n$ is never large, because of the
exponentially growing factor $2^n$ in the complexity, so the difference
between the two approaches is essentially an implementation-dependent constant
factor.

\section{Parallel algorithms}

Suppose we wish to use $P \ge 1$ processors in parallel. If the
$K \simeq 2^n/n$ necklaces of length $n$ are 
$W_0 = 0\ldots 0, W_1, W_2, \ldots, W_{K-1} = 1\ldots 1$,
we would like to ask processor $q$ ($0 \le q < P$) to compute the determinants
corresponding to necklaces
$W_{\lfloor{qK/P}\rfloor},\ldots,W_{\lfloor{(q+1)K/P}\rfloor-1}$.
The problem is how to determine
the starting point $W_{\lfloor{qK/P}\rfloor}$ for processor~$q$,
without enumerating $W_1, W_2, \ldots, W_{\lfloor{qK/P}\rfloor}$.
A polynomial-time algorithm for this problem is claimed in~\cite{Kociumaka14},
but it is very complicated. We preferred to adopt a simpler approach
which is much easier to implement and sufficient in practice.%
\footnote{A similar algorithm, although not specifically intended for
parallel computation, is described in~\cite{Briggs18}.
For an algorithm using less storage, and related references,
see~\cite{Ghosh13}.}

The idea is to take a random sample of (say) $T := 4000P^2$ necklaces
(each of length $n$). Sort the sample, and then divide it into $P$
equal-sized segments. Modify the initial segment to start with 
$W_0 = 0\ldots 0$
and the final segment to end with $W_{K-1} = 1\ldots 1$.  Thus, each
processor has the same number $\lfloor K/P \rfloor$ words to process, apart
from a small sampling error which is negligible in practice.
Also, we know
the necklace starting each segment, so we can use Duval's algorithm to
enumerate all necklaces in a segment.

We describe how to randomly sample the set of all
necklaces of length $n$ in such a manner that each necklace occurs in
the sample with equal probability.  Generate a random binary string
of length $n$, and test (using Booth's algorithm) if it corresponds to a
necklace. If so, the string is accepted.
Otherwise, the string is rejected and we try again.  The process
is repeated until we have the desired number $T$ of necklaces (not
necessarily distinct).  Clearly each
necklace is equally likely to appear in the final list.
Since the probability that a random binary string is a necklace
is close to $1/n$, the number of random binary strings that are needed
is of order $nT$.  What we have described is a simple example of
Von Neumann's \emph{rejection method}, first described by Forsythe 
in~\cite{von-Neumann-work}.
Other examples may be found in Devroye's book~\cite{Devroye86}.

\section{Computational results}

In Tables \ref{tab:a}--\ref{tab:b} 
we give computational results for the 
maximal determinants $D_{01}(n)$ of
$\{0,1\}$-circulants of order $n \le 53$.  The third column of each table
gives the ratio $D_{01}(n)/U_{01}(n)$, where $D_{01}(n)$ is the maximum
of $|\det(B)|$ for $\{0,1\}$-circulants $B$ of order $n$, and
$U_{01}(n)$ is an upper bound on $D_{01}(n)$.

Similarly, in Tables \ref{tab:c}--\ref{tab:d} 
we give computational results for the maximal
determinants $D_{\pm1}(n)$ of $\{\pm1\}$-circulants of order $n \le 53$.
Here the third column is the ratio $D_{\pm1}(n)/U_{\pm1}(n)$,
where $U_{\pm1}(n)$ is an upper bound on $D_{\pm1}(n)$.
In Tables \ref{tab:c}--\ref{tab:d}
we scale the determinants of $\{\pm1\}$-circulants by
dividing by the known factor $2^{n-1}$.
In the last column of Table~\ref{tab:c}, ``$-$'' and ``$+$'' are used
as abbreviations for $-1$ and $+1$ respectively. 

The bounds $U_{01}(n)$ and $U_{\pm1}(n)$ are defined as follows.  Let
\begin{equation}
\HBE(n) := \begin{cases}
	    n^{n/2} \text{ if } n \equiv 0 \bmod 4,	\\
	    2(n-1)\,(n-2)^{(n-2)/2} \text { if } n \equiv 2 \bmod 4,	\\
	    (2n-1)^{1/2}\,(n-1)^{(n-1)/2} \text{ otherwise.}
	    \end{cases}
							\label{eq:HBE}
\end{equation}
Then $\HBE(n)$ is an upper bound on $|\det A|$ for $\{\pm1\}$-matrices $A$
of order~$n$. The case $n \equiv 0 \bmod 4$ is due to
Hadamard~\cite{Hadamard93};
the case $n \equiv 2 \bmod 4$ is due to 
Ehlich~\cite{Ehlich64a} and Wojtas~\cite{Wojtas64};
and the remaining case ($n$ odd)
is due to Barba~\cite{Barba33}, Ehlich~\cite{Ehlich64a},
and Wojtas~\cite{Wojtas64}.
We do not use Ehlich's slightly sharper, but more complicated, bound
that applies when $n \equiv 3 \bmod 4$.
For this bound, see Ehlich~\cite{Ehlich64b} or Orrick~\cite{Orrick-www}.

In view of the discussion in \S\ref{sec:intro}, we take
\begin{equation}					\label{eq:pmbound}
U_{\pm1}(n) := 2^{n-1}\lfloor{\HBE(n)/2^{n-1}}\rfloor
\end{equation}
and
\begin{equation}					\label{eq:01bound}
U_{01}(n) := \lfloor\HBE(n+1)/2^n \rfloor.
\end{equation}

It is an open question whether $D_{\pm1}(n)$ attains the
bound $U_{\pm1}(n)$ for any $n > 13$.  (If we restrict attention to
the cases $n \equiv 0 \bmod 4$, this is the \emph{circulant Hadamard}
problem.)  On the other hand, $D_{01}(p) = U_{01}(p)$ for all primes
$p \equiv 3 \bmod 4$.  This follows from the first \emph{Paley
construction}~\cite{Paley33}, 
which constructs a Hadamard matrix of order $p+1$
with a circulant submatrix of order $p$. Inspection of
Tables \ref{tab:a}--\ref{tab:b} reveals
that $D_{01}(n) = U_{01}(n)$ in some other cases, specifically
$n \in \{1, 2, 4, 15, 35\}$.

Table \ref{tab:b} extends the list of $D_{01}(n)$ values given for
$n \le 37$ in OEIS A086432 and the associated b-file~\cite{A086432}.
Table \ref{tab:d} extends the list of $D_{\pm1}(n)/2^{n-1}$ values given for
$n \le 28$ in OEIS A215897~\cite{A215897}. This implies a 
corresponding extension for OEIS A215723, which lists the
unscaled values $D_{\pm1}(n)$.

As an indication of the time required to compute the tables, we note
that the computation of $D_{01}(52)$ using our parallel program
(implemented in C using GMP~\cite{GMP}) took $11$ processor-years
using 128 Intel Xeon3 (2.2GHz) and 224 Xeon4 (2.6GHz) processors. 
The computation time for order $n$ was roughly proportional to $2^n$. 

For verification, all the values given in the
tables for orders $n < 50$
were computed at least twice, using different programs and/or
different prime moduli~$p$.
No discrepancies were found.

\pagebreak[3]

\section{Some conjectures}

In this section we discuss, and disprove, some plausible conjectures.

\subsection*{Conjecture A}

{From} the third column of Table~\ref{tab:a}, the determinant
of a $\{0,1\}$-circulant can attain the upper bound $U_{01}(n)$ 
in the cases \hbox{$n \in \{1,2,3,4,7,11,15,19,23\}$}. The Paley construction
explains this for $n = 3, 7, 11, 19, 23$, and larger cases where $n$
is a prime and $n \equiv 3 \bmod 4$.  However, it does not explain
the case $n=15=3 \times 5$. 
Also, the upper bound is not attained for $n = 27 = 3^3$.   
Thus, a plausible conjecture is that the upper bound can be attained
whenever $n \equiv 3 \bmod 4$ is 
the product of at most two distinct primes. Support is provided by the
computation for $n = 35 = 5 \times 7$, since
$D_{01}(35) = U_{01}(35)$.

Our computation for $n=39$ disproves these conjectures,
since $39 = 3 \times 13$ is a product of two distinct primes, but
$D_{01}(39) < U_{01}(39)$.  Another counter-example is $n=51 = 3 \times 17$.
We find that $D_{01}(51) < U_{01}(51)/2$.

After our computations were completed, we discovered an explanation
for the cases $n=15$ and $n=35$.  In each case $n$ has the form
$p(p+2)$, where $p$ and $p+2$ are both prime. Such $n$ are covered
by case (2) 
of the following theorem, which we quote (with a slight change in notation)
from~\cite{KKS06}.
Note that a ``circulant core'' of order $n$ refers to a
$\{0,1\}$-circulant matrix of order $n$ which can be used to construct
a Hadamard matrix of order $n+1$ using the correspondence between
$\{0,1\}$-matrices and $\{\pm1\}$-matrices described in~\S1.
\begin{theorem}[Hadamard circulant core construction]
\label{thm:circulant_core}
A Hadamard matrix of order $n+1$ with circulant core of order $n$ 
exists if
\begin{itemize}
\item[(1)] $n \equiv 3 \bmod 4$ is a prime;
\item[(2)] $n = p(p+2)$, where $p$ and $p+2$ are prime;
\item[(3)] $n = 2^k-1$, where $k$ is a positive integer; or
\item[(4)] $n = 4k^2 + 27$, where $k$ is a positive integer
		and $n$ is a prime.
\end{itemize}
\end{theorem}
\begin{proof}
Case (1) is due to Paley~\cite{Paley33}; 
case (2) is due to Stanton and Sprott~\cite{Stanton58} and
also Whiteman~\cite{Whiteman62};
case (3) is due to Singer~\cite{Singer38};
and case (4) is due to Hall~\cite[Theorem 2.2]{Hall56}.
\end{proof}

Hall~\cite[pg.~$980$]{Hall56} remarks
that case~(4) 
is subsumed by case~(1), since $4k^2+27 \equiv 3 \bmod 4$, 
but we mention case~(4) since
Hall's construction is different from that of Paley.

We do not know if the list given by Theorem~\ref{thm:circulant_core} is
exhaustive.  The computational results given in
Tables~\ref{tab:a}--\ref{tab:b} show that, for $1 \le n \le 53$,
only those $n$
given by Theorem~\ref{thm:circulant_core} can provide a Hadamard matrix of
order $n+1$ with a circulant core.
Also, a circulant $\{0,1\}$-matrix of order $n \le 53$ can achieve the upper
bound \eqref{eq:01bound} if and only if $n \le 4$ or $n$ satisfies condition
(1), (2) or (3) of Theorem~\ref{thm:circulant_core}.
 
\subsection*{Conjecture B, case $[0,1]$}

When considering maximal determinants of matrices with real elements
in the interval $[0,1]$, we can see that the maximum occurs at
extreme points of the polytope.\footnote{This is already implicit in
Hadamard~\cite{Hadamard93}.}
To prove this, we need only note that the determinant $\det A$
of a square matrix $A= (a_{j,k})$ is a linear function of each
variable $a_{j,k}$ considered separately.  Thus, if a local maximum
of $\det A$ occurs for some $a_{j,k} \in (0,1)$, we can
replace $a_{j,k}$ by (at least one of) $0$ or $1$ without decreasing
$\det A$.

This argument does not apply if $A$ is restricted to be a circulant of
order $n > 1$, because
then the free parameters are just the elements $a_0,\ldots,a_{n-1}$ of the
first row of $A$, and $\det A$ is \emph{not} a linear function of
each $a_j$ considered separately. For example, if $n=2$ we have
$\det A = a_0^2 - a_1^2$. Nevertheless, inspection of small cases suggests
the conjecture that the maximum of $|\det A|$ occurs at extreme points of
the $n$-dimensional polytope.

We were unable to prove the conjecture, so wrote a program to check it
numerically, and found that, in general, the conjecture is false.

The idea is as follows.  Consider all possible circulants $A$
of order $n$ with entries in $\{0,1\}$. 
If $\det A = \pm D_{01}(n)$, check if a small perturbation
of $a_0$ towards the interior
of the polytope would increase $|\det A|$.  Although such behaviour is
rare, it does occur.\footnote{For reasons of efficiency, 
our program takes as input a list (generated during the computation of
Tables \ref{tab:a}--\ref{tab:b}) of necklaces
that define circulants $A$ with maximal $|\det A|$, then considers all
possible rotations of these circulants.}

The smallest examples occur for $n=9$.
Consider $A = \circulant(a_0,\ldots, a_8)$ with
$(a_0,\ldots,a_8) = (0,0,0,1,1,1,1,0,1)$. We have
$\det A = 95 = D_{01}(9)$, but
$\partial \det A/\partial a_0 = 9$.
If $a_0 = \x$ for some small $\x$,
then $|\det A(\x)| = 95 + 9\x + O(\x^2)$,
so $|\det A(\x)| > 95$ for sufficiently small
$\x > 0$.
In fact, 
$|\det A(0.241)| > 96.757$.

For $n=10$, an example is
$A = \circulant(0, 0, 1, 0, 0, 1, 1, 1, 1, 0)$, $\det A = 275$.
Replacing $a_0$ by $\x = 0.112$, we obtain $\det A(\x) > 279.4$.

We found examples of such behaviour for $n=9, 10$ and no other $n$
up to the limit of Table~\ref{tab:b}.
However, there is a different class of examples that occur when 
$n = 4k+1 > 5$ is a prime, e.g.~$n = 13, 17, 29, 37, 41$, etc.
For this class we make a small modification to the \emph{Uniformly Redundant
Arrays} (URAs) of \cite{Busboom97,Fenimore78}, which are equivalent to
Abelian difference sets~\cite{Hall56}.\footnote{Our construction
is also close to the ``modified'' URAs (MURAs)
of~\cite{Gottesman89}.}
Define
\[
A_n(x) := \circulant\left(x,
		\textstyle\frac{1+\chi(1)}{2},\frac{1+\chi(2)}{2},\ldots,
		\frac{1+\chi(n-1)}{2}\right), 
\]
where $\chi$ is a quadratic character, defined by the Legendre symbol
\[
\chi(j) = \left(\frac{j}{n}\right) := 
 \begin{cases}
 +1 \text{ if $j$ is a quadratic residue modulo $n$ and $j \ne 0 \bmod n$};\\
 -1 \text{ if $j$ is a quadratic non-residue modulo $n$};\\
 \phantom{+}0 \text{ if $j \equiv 0 \bmod n$}.
 \end{cases}
\]
Then $A_n(0)$ corresponds to a $1$-D URA, but $\det A_n(0)$ is not generally
maximal in the class of circulant determinants.
However, $\det A_n(\frac12)$ may be larger
than the corresponding entry in Tables~\ref{tab:a}--\ref{tab:b}.
It may be shown\footnote{The proof uses the identity
$A_n(0)^2 + A_n(0) = k(I+J)$.}
that, for $n = 4k+1$ an odd prime,
\begin{equation}
\det A_n(x) = (x+2k)(x^2-x-k)^{2k}.		\label{eq:detAx}
\end{equation}
In particular,
$\det A_n(0) = 2 k^{2k+1}$,
$\det A_n(1) = (2k+1)k^{2k}$,
and
\[
\det A_n({\textstyle\frac12}) = 2^{-n}\,n^{(n+1)/2}.
\]
It may be verified numerically that
$\det A_n(\frac12)$ exceeds the maximal determinant
given in Tables~\ref{tab:a}--\ref{tab:b} for $n = 13, 17, 29, 37,$ and $41$.
The next possibility, $n=53$, is beyond the range of Table~\ref{tab:b}.

We observe that the maximum of $\det A_n(x)$ for $x\in [0,1]$
is not at $x=\frac12$.
One can show, by logarithmic differentiation of~\eqref{eq:detAx},
that a local maximum occurs at
\[
x = x_k := \frac{\sqrt{1+4k^2}+1-2k}{2} = \frac12 + \frac{1}{8k}
	+ O(k^{-3}),
\]
and
\[
 \max_{0\le x \le 1}\det A_n(x)
 = \det A_n(x_k)
 = \det A_n({\textstyle\frac12})\left(1+\frac{1}{8kn} + O(k^{-4})\right).
\]
For example, if $k=3, n=13$, we have 
$x_3 = (\sqrt{37}-5)/2 \approx 0.5414$,
and 
$	U_{01}(13) = 9477
	> \det A_{13}(x_3) \approx 7684.16
	> \det A_{13}({\textstyle\frac12}) \approx 7659.73
	> D_{01}(13) = 6561
	> \det A_{13}(1) = 5103
	> \det A_{13}(0) = 4374.
$

\subsection*{Conjecture B, case $[-1,1]$}

Replacing $[0,1]$ by $[-1,1]$, we find analogous behaviour for
$n = 2, 9, 10, 11, 18$, $22$ and no other $n$ up to the limit
of Table~\ref{tab:d}.
The case $n=2$ is trivial because,
for circulants of order~$2$ over $\Spm$, we necessarily have
$\det A = 0$ at the extreme points
$(a_0,a_1) = (\pm1, \pm1)$.

The other cases are non-trivial.
For example, if $n=9$, consider
\[A(\x) := \circulant(1-\x,1,-1,1,-1,-1,1,1,1).\]
We find that
\[\det A(\x) = 6912 + 4608\x + O(\x^2),\]
so sufficiently small $\x > 0$ gives $\det A(\x) > 6912 = D_{\pm1}(9)$.
Indeed, we can take $\x=1$, as $\det A(1) = 8582 > 6912$.

If $n=10$, we find that
\[\det\circulant(1-\x, -1, 1, 1, -1, -1, -1, -1, -1, -1)
 = -(22528 + 2560\x + O(\x^2)),\]
and
\[\det\circulant(-1+\x, -1, -1, 1, -1, 1, 1, -1, -1, -1)
 = 22528 + 7680\x + O(\x^2),\]
so in both cases a sufficiently small $\x>0$ disproves the conjecture.
A different type of exceptional case is illustrated by
\[
A(x) := \circulant(x, -1, 1, -1, 1, 1, -1, -1, -1, -1),
\]
where we find that $\det A(x)$ is an even polynomial in $x$,
and 
\[-\det A(0) = 33489 > -\det A(\pm1) = 22528 = D_{\pm1}(10).\]

Similarly, for order $22$, consider
\[
A(x) := \circulant(x,\! -1,\! 1,\! 1,\! -1,\! -1,\! -1,\! -1,\! -1,\! -1,\! -1,\! 1,\! 1,\! 
	-1,\! 1,\! -1,\! 1,\! -1,\! 1,\! 1,\! -1,\! -1).
\]
Then 
\[
-\det A(0) = 216409254831025 > -\det A(\pm1) = 215055782117376.
\]
Since $215055782117376 = D_{\pm1}(22) = 2^{21}\times 102546588$ (see
Table~\ref{tab:c}), we have $|\det A(0)| > D_{\pm1}(22)$.

Our search was not exhaustive, so there
may be other $n$ within the range of Tables~\ref{tab:c}--\ref{tab:d} 
for which the maximum determinant does not occur at
an extreme point of $[-1,1]^n$.

\section{Remarks on periodic autocorrelations}

It is hard to discern a pattern in the lex-least words given in
Tables~\ref{tab:a}--\ref{tab:d}.  It seems more fruitful to consider
the \emph{periodic autocorrelations} of the first rows of the corresponding
circulants. Equivalently, we can consider the first rows of the
\emph{Gram matrices} $G = A^TA$, where $A$ is the relevant circulant.

In the case of $(0,1)$-circulants, it may be useful to map 
$(0,1) \mapsto (-1,1)$, and consider the first row of 
$G' = (2A-J)^T (2A-J)$.
For example, provided $n > 4$, 
the upper bound is achieved in Tables \ref{tab:a}--\ref{tab:b}
if and only if the first row of $G'$ is 
$(n,-1,-1,\ldots,-1)$, see~\cite{MacWilliams76}.

In some cases the maximal determinants given in 
Tables~\ref{tab:a}--\ref{tab:d} have only small prime factors.
For example, the entry for $n=52$ in Table~\ref{tab:d} is
$2^{49}\, 3^{24}\, 5^4$, and this can be explained if we observe that
the first row of $G$ is 
$(52, 0, 0, 0, 4, 0, 0, 0, 4, \ldots, 0, 0, 0, 4, 0, 0, 0)$.
Thus, we can write $G = 52I + 4E^4 + 4E^8 + \cdots + 4E^{48}$, where $E$
is the ``circular shift'' matrix. 
Similarly, the entry for $n=48$ in Table~\ref{tab:d} is
$2^{49}\,3^6\,5^{12}$,
and here $G = 48I + 4E^{12} + 8E^{24} + 4E^{36}$.

\section{Acknowledgements}

We thank J\"org Arndt for his comments on a draft of this document, Alex Arkhipov 
for his helpful comments about Hadamard matrices with circulant cores and related
group theory,
and the authors of Magma~\cite{Magma1} and GMP~\cite{GMP} for their
excellent software. Computing resources were provided by the Australian
National University and the University of Newcastle (Australia). The first
author was supported in part by Australian Research Council grant
DP140101417.

\pagebreak[3]

\bigskip
\hrule
\bigskip

\pagebreak[3]

\noindent 2010 {\it Mathematics Subject Classification}:
Primary 05A15;
Secondary 05A19, 65T50.

\noindent \emph{Keywords: } 
binary matrix,
Booth's algorithm,
circulant,
circulant core,
computational imaging,
convolutional Gaussian channels,
difference set,
discrete Mahler measure,
Duval's algorithm,
Hadamard bound,
Hadamard matrix,
Lyndon word,
maximal determinant,
modular computation,
MURA,
necklace,
parallel algorithm,
parallel computation,
URA

\bigskip
\hrule
\bigskip

\noindent (Concerned with sequences
\seqnum{A000031},
\seqnum{A086432},	
\seqnum{A215723},	
\seqnum{A215897}.)	

\bigskip
\hrule
\bigskip

\vspace*{+.1in}
\noindent

\noindent

\pagebreak[4]

\section*{Appendix -- Tables of Maximal Determinants}

\begin{table}[ht] 	
\vspace*{-0pt}
\begin{center}
\begin{tabular}{c|c|c|cc}
order	       & maximal 	& ratio  & lex-least & lex-least \\ 
$n$	  &  $|$determinant$|$	& $D_{01}(n)/$ & word & word \\
 	  & $D_{01}(n)$ & $U_{01}(n)$ & (decimal) & (over $\{0,1\}$) \\
\hline
$1$ & $1$ & $1.0000$ & $1$ & \nsp{\tt 1} \\	
$2$ & $1$ & $1.0000$ & $1$ & \nsp{\tt 01} \\
$3$ & $2$ & $1.0000$ & $3$ & \nsp{\tt 011} \\
$4$ & $3$ & $1.0000$ & $7$ & \nsp{\tt 0111} \\
$5$ & $4$ & $0.8000$ & $15$ & \nsp{\tt 01111} \\
$6$ & $9$ & $0.7500$ & $11$ & \nsp{\tt 001011} \\
$7$ & $32$ & $1.0000$ & $23$ & \nsp{\tt 0010111} \\
$8$ & $45$ & $0.6923$ & $47$ & \nsp{\tt 00101111} \\
$9$ & $95$ & $0.6597$ & $47$ & \nsp{\tt 000101111} \\
$10$ & $275$ & $0.6152$ & $55$ & \nsp{\tt 0000110111} \\
$11$ & $1458$ & $1.0000$ & $183$ & \nsp{\tt 00010110111} \\
$12$ & $2240$ & $0.6145$ & $439$ & \nsp{\tt 000110110111} \\
$13$ & $6561$ & $0.6923$ & $1527$ & \nsp{\tt 0010111110111} \\
$14$ & $19952$ & $0.5759$ & $751$ & \nsp{\tt 00001011101111} \\
$15$ & $131072$ & $1.0000$ & $2479$ & \nsp{\tt 000100110101111} \\
$16$ & $214245$ & $0.5691$ & $2935$ & \nsp{\tt 0000101101110111} \\
$17$ & $755829$ & $0.6784$ & $2935$ & \nsp{\tt 00000101101110111} \\
$18$ & $2994003$ & $0.6505$ & $9903$ & \nsp{\tt 000010011010101111} \\
$19$ & $19531250$ & $1.0000$ & $22427$ & \nsp{\tt 0000101011110011011} \\
$20$ & $37579575$ & $0.6010$ & $28023$ & \nsp{\tt 00000110110101110111} \\
$21$ & $134534444$ & $0.6560$ & $45999$ & \nsp{\tt 000001011001110101111} \\
$22$ & $577397064$ & $0.6178$ & $117623$ & \nsp{\tt 0000011100101101110111} \\  
$23$ & $4353564672$ & $1.0000$ & $340831$ & \nsp{\tt 00001010011001101011111} \\
$24$ & $10757577600$ & $0.7060$ & $843119$ & \nsp{\tt 000011001101110101101111} \\
$25$ & $31495183733$ & $0.5787$ & $638287$ & \nsp{\tt 0000010011011110101001111} \\
\hline
\end{tabular}
\end{center}
\caption{Maximal determinants of $\{0,1\}$-circulants of order $n \le 25$.}
\vspace*{-20pt}
\label{tab:a}
\end{table}

\begin{table}[ht]
\begin{center}
\begin{tabular}{c|c|c|c}
order     & maximal  & ratio to & lex-least word \\
$n$	& $|$determinant$|$	& upper bound & (decimal) \\
\hline
$26$ & $154611524732$ & $0.5744$ & $957175$ \\
$27$ & $738139162166$ & $0.5442$ & $1796839$ \\
$28$ & $3124126889325$ & $0.6101$ & $5469423$ \\
$29$ & $11937232425585$ & $0.6069$ & $6774063$ \\
$30$ & $65455857159975$ & $0.6271$ & $37463883$ \\
$31$ & $562949953421312$ & $1.0000$ & $77446231$ \\
$32$ & $1395230053365015$ & $0.6148$ & $47828907$ \\
$33$ & $5687258414265018$ & $0.6123$ & $196303815$ \\
$34$ & $30551195956571643$ & $0.5827$ & $95151003$ \\
$35$ & $300189270593998242$ & $1.0000$ & $1324935477$ \\
$36$ & $809028975189744400$ & $0.6309$ & $1822895095$ \\
$37$ & $3198686446402685263$ & $0.5760$ & $430812063$ \\
$38$ & $19288701806345611347$ & $0.5825$ & $2846677239$ \\
$39$ & $103227456252120723684$ & $0.5161$ & $10313700815$ \\
$40$ & $529663503370085366373$ & $0.5885$ & $6269629671$ \\
$41$ & $2311393009109010944326$ & $0.5638$ & $26764629467$ \\
$42$ & $15469925980869995489631$ & $0.6023$ & $22992859983$ \\
$43$ & $162805498773679522226642$ & $1.0000$ & $92035379515$ \\
$44$ & $402826140168935435652453$ & $0.5245$ & $162368181483$ \\
$45$ & $2268175963362305735661143$ & $0.6192$ & $226394696439$ \\
$46$ & $12738408112895861486972391$ & $0.5307$ & $631304341299$ \\
$47$ & $158993694406781688266883072$ & $1.0000$ & $4626135339999$ \\
$48$ & $483776963047101724429782080$ & $0.6179$ & $924925407055$ \\
$49$ & $2226275734022433928055705600$ & $0.5715$ & $1588449170843$ \\
$50$ & $15940963431893953997118039375$ & $0.5992$ & $5455102172067$ \\
$51$ & $86343902346653136953496818019$ & $0.4706$ & $12463552538547$ \\
$52$ & $471252255596620483490068604560$ & $0.5013$ & $23418838481755$ \\
$53$ & $2670231923706326010918104225583$ & $0.5492$ & $12803059922743$ \\ 
\hline
\end{tabular}
\end{center}
\caption{Maximal determinants of $\{0,1\}$-circulants, $25 <  n \le 53$.}
\vspace*{-10pt}
\label{tab:b}
\end{table}

\begin{table}[ht] 	
\begin{center}
\begin{tabular}{c|c|c|cc}
order	       & maximal  	& ratio  & lex-least & lex-least \\ 
$n$	  &  scaled $|$det$|$	& $D_{\pm}(n)/$ & word & word \\
 	  & $D_{\pm1}(n)/2^{n-1}$ & $U_{\pm}(n)$ & (decimal) & (over $\{\pm1\}$) \\
\hline
1 & 1 & 1.0000 & 0 & \nsp{\tt -} \\
2 & 0 & 0.0000 & 0 & \nsp{\tt --} \\
3 & 1 & 1.0000 & 1 & \nsp{\tt --+} \\
4 & 2 & 1.0000 & 1 & \nsp{\tt ---+} \\
5 & 3 & 1.0000 & 1 & \nsp{\tt ----+} \\
6 & 4 & 0.8000 & 1 & \nsp{\tt -----+} \\
7 & 8 & 0.6667 & 11 & \nsp{\tt ---+-++} \\
8 & 18 & 0.5625 & 11 & \nsp{\tt ----+-++} \\
9 & 27 & 0.4154 & 11 & \nsp{\tt -----+-++} \\
10 & 44 & 0.3056 & 11 & \nsp{\tt ------+-++} \\
11 & 267 & 0.5973 & 39 & \nsp{\tt -----+--+++} \\
12 & 1024 & 0.7023 & 83 & \nsp{\tt -----+-+--++} \\
13 & 3645 & 1.0000 & 83 & \nsp{\tt ------+-+--++} \\
14 & 6144 & 0.6483 & 83 & \nsp{\tt -------+-+--++} \\
15 & 23859 & 0.6886 & 359 & \nsp{\tt ------+-++--+++} \\
16 & 50176 & 0.3828 & 691 & \nsp{\tt ------+-+-++--++} \\
17 & 187377 & 0.4977 & 1643 & \nsp{\tt ------++--++-+-++} \\
18 & 531468 & 0.4770 & 2215 & \nsp{\tt ------+---+-+--+++} \\
19 & 3302697 & 0.7176 & 9895 & \nsp{\tt -----+--++-+-+--+++} \\
20 & 10616832 & 0.5436 & 6483 & \nsp{\tt -------++--+-+-+--++} \\
21 & 39337984 & 0.6291 & 67863 & \nsp{\tt ----+----+--+---+-+++} \\
22 & 102546588 & 0.5000 & 21095 & \nsp{\tt -------+-+--+--++--+++} \\ 
23 & 568833245 & 0.6087 & 72519 & \nsp{\tt ------+---++-++-+---+++}\\
24 & 3073593600 & 0.7060 & 144791 & \nsp{\tt ------+---++-+-++--+-+++} \\
25 & 8721488875 & 0.5724 & 108199 & \nsp{\tt --------++-+--++-+-+--+++}\\
\hline
\end{tabular}
\end{center}
\caption{Maximal scaled determinants of $\{\pm1\}$-circulants
 of order $n \le 25$.}
\vspace*{-10pt}
\label{tab:c}
\end{table}

\begin{table}[ht]
\begin{center}
\begin{tabular}{c|c|c|c}
order     & maximal scaled $|$det$|$ 	& ratio to & lex-least word\\
 $n$	  &  $D_{\pm1}(n)/2^{n-1}$	& upper bound & (decimal)\\
\hline
26 & 32998447572 & 0.6064 & 355463 \\
27 & 164855413835 & 0.6125 & 604381 \\
28 & 572108938470 & 0.4218 & 1289739 \\
29 & 2490252810073 & 0.4863 & 1611219 \\
30 & 10831449635712 & 0.5507 & 1680711 \\
31 & 68045615234375 & 0.6520 & 6870231 \\
32 & 282773291271138 & 0.5023 & 12817083 \\
33 & 1592413932070703 & 0.7017 & 18635419 \\
34 & 5234078743146888 & 0.5635 & 55100887 \\
35 & 33374247484277975 & 0.6366 & 149009085 \\
36 & 198124573871046186 & 0.6600 & 160340631 \\
37 & 787413957917252603 & 0.6140 & 415804239 \\
38 & 3195257068570067448 & 0.5754 & 829121815 \\
39 & 22999238901574021485 & 0.6946 & 4737823097 \\
40 & 117140061677844350646 & 0.5857 & 1446278811 \\
41 & 536469708946538168543 & 0.5961 & 3001209959 \\
42 & 2417648227367853639168 & 0.5897 & 19153917469 \\
43 & 14611334654738350617599 & 0.5689 & 52222437727 \\
44 & 65738632907943707712320 & 0.4038 & 20159598251 \\
45 & 438910341492340511320163 & 0.5715 & 166482220965 \\
46 & 2010768410464246499566152 & 0.5489 & 90422521191 \\
47 & 12779930756727248097293989 & 0.5324 & 115099593371 \\
48 & 100192997081088000000000000 & 0.6302 & 242235026743 \\
49 & 408375323859124630659059549 & 0.5216 & 1416138805685 \\
50 & 2152519997519833685106486024 & 0.5526 & 2380679727935 \\
51 & 14098690136202107270366810369 & 0.5300 & 2716242515341 \\
52 & 99371059004238555166801920000 & 0.5416 & 1758408815375 \\
53 & 512364770126478307153560491081 & 0.5451 & 10146024354919 \\
\hline
\end{tabular}
\end{center}
\caption{Maximal scaled determinants of $\{\pm1\}$-circulants,
 $25 <  n \le 53$.}
\vspace*{-10pt}
\label{tab:d}
\end{table}

\end{document}